%% file: main.tex
\documentclass{llncs}

\usepackage{amsfonts}
\usepackage{amsmath}
\usepackage{booktabs} 
\usepackage{tikz}
\usetikzlibrary{calc}
\usetikzlibrary{matrix}
\usepackage{xcolor}
\usepackage{graphicx}
\usepackage{mathtools}

\usepackage{bm}
\usepackage{algorithm}
\usepackage[noend]{algpseudocode}

\newcommand{\red}[1]{\textcolor{red}{#1}}

\usepackage{soul}
\usepackage{etoolbox}

\definecolor{webred}{rgb}{0.5,0,0}
\definecolor{webblue}{rgb}{0,0,0.8}
\usepackage[pdftex,colorlinks,citecolor=webblue,linkcolor=webred,]{hyperref}

\usepackage[caption=false]{subfig}

\newcommand{\norm}[1]{\left\lVert #1 \right\rVert}

\newcommand{\new}[1]{#1}
\newcommand{\old}[1]{}

\title{Numerical Verification of Affine Systems \\with up to a Billion Dimensions}

\author{Stanley Bak \inst{1}, Hoang-Dung Tran \inst{2}, and Taylor T. Johnson \inst{2}}
\institute{Safe Sky Analytics \and Vanderbilt University}

\begin{document}

\maketitle

\begin{abstract}
\input{abstract}
\end{abstract}

\makeatletter
\def\blfootnote{\xdef\@thefnmark{}\@footnotetext}
\makeatother
\blfootnote{DISTRIBUTION A. Approved for public release; Distribution 
unlimited. \\AFRL PA \# 88ABW-2017-5562 cleared on 07 Nov 2017.\vspace{-2em}}

\input{intro}
\input{basic}

\input{memory}
\input{computational}
\input{scalability}

\input{evaluation}

\input{related_work}
\input{conclusion}

\bibliographystyle{abbrv}
\bibliography{bak,tran}


\end{document}

%% file: abstract.tex
Affine systems reachability is the basis of many verification methods.
With further computation, methods exist to reason about richer models with inputs, 
nonlinear differential equations, and hybrid dynamics.
As such, the scalability of affine systems verification is a prerequisite to scalable analysis for more complex systems.
In this paper, we improve the scalability of affine systems verification, in terms of
the number of dimensions (variables) in the system.
\vspace{1.0em}

The reachable states of affine systems can be written in terms of the matrix exponential, and safety checking can be performed at specific time steps with linear programming.
Unfortunately, for large systems with many state variables, this direct approach requires an intractable amount of memory while using an intractable amount of computation time.
We overcome these challenges by combining several methods that leverage common problem structure.
Memory is reduced by exploiting initial states that are not full-dimensional and safety properties (outputs) over a few linear projections of the state variables.
Computation time is saved by using numerical simulations to compute only projections of the matrix 
exponential relevant for the verification problem.
Since large systems often have sparse dynamics, we use Krylov-subspace simulation approaches based on
the Arnoldi or Lanczos iterations.
Our method produces accurate counter-examples when properties are violated and, in the extreme case with sufficient problem structure, can analyze a system with one billion real-valued 
state variables.

%% file: intro.tex
\section{Introduction}
\label{sec:intro}

An affine system is modeled with the ordinary differential equation $\dot{x} = Ax + b$,
where $x$ is a vector of $n$ state variables, $A$ is the
$n \times n$ dynamics matrix, and $b$ is an $n \times 1$ vector of constant forcing terms.
Given a set of initial states, a set of unsafe states, and a time bound,
the time-bounded safety verification problem is to check if there exists an initial state and a time
within the bound such that the solution of the affine system enters the unsafe set.
%
%

One way to solve the verification problem is to construct the reachable set of states.
The reachable set contains all states that lie along any solution to the differential equation $\dot{x} = Ax + b$,
starting from any initial state up to the time bound.
If the reachable set does not intersect with the unsafe states, then the system is verified as safe.
In the discrete-time setting, we construct the reachable set at each time instant, and then
perform the unsafe check once per step using linear programming (LP).

This discrete-time approach forms the core of many verification methods for
richer classes of systems.
For example, when systems have time-varying inputs, an additional input-effects term can be
computed at each step and added to the discrete-time reachable set using a Minkowski sum
operation~\cite{girard2006efficient,bak2017cav}.
Overapproximation of the continuous-time reachable set
is possible by noting that, in the finite time
between time steps, the system can only go a bounded distance from the discrete-time solution.
Based on this observation, methods exist that perform bloating from the discrete-time solution in order to guarantee an
overapproximation of the continuous-time reachable set~\cite{spaceex,le2010reachability,le2009reachability}.
The reachable set for nonlinear dynamical systems can also be overapproximated with techniques based on
affine methods,
by linearizing the nonlinear dynamics and then adding uncertain terms to account for mismatch between the
linear and nonlinear systems.
In hybridization methods~\cite{dang10,bak2016hscc,althoff2008reachability}, this process is repeated in multiple
domains to reduce the overapproximation error.
Finally, methods to verify hybrid systems that combine continuous dynamics
and discrete behaviors, such as a physical system controlled by software, also build upon
the core operations needed to analyze affine systems~\cite{maler1991hybrid,alur1993hybrid},
in conjunction with additional techniques to handle combinatorial aspects.
%
%
All of these powerful methods build on the core computations used for affine systems reachability.
In this paper, we focus on the scalability of this fundamental computation.

Verification approaches for systems that have real numbers can be categorized into \emph{validated} methods
and \emph{numerical} methods.
Validated methods, such as interval analysis~\cite{stauning1997automatic},
maintain guaranteed bounds on values used throughout the computation.
Numerical methods, on the other hand, accept using finite-precision floating-point numbers and algorithms that
perform operations up to any user-desired accuracy, such as finite series expansions to compute a matrix exponential.
Although desirable, validated methods are typically slower and often fail to work on large systems to due
the accumulation of overapproximation error.
In this work, we focus on numerical verification methods, as the scale of systems we want to analyze would make current validated
approaches infeasible.

\begin{figure}[t]
    \centering
    \includegraphics[width=0.7\columnwidth]{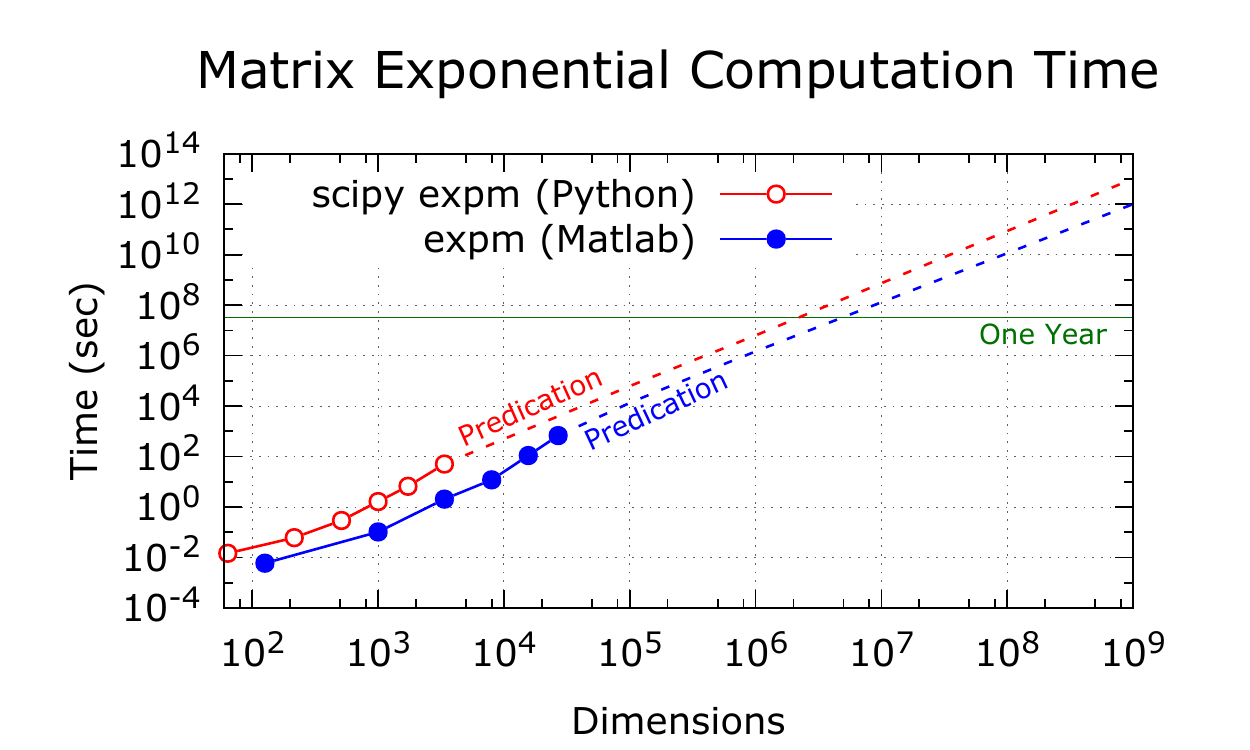}%
		\vspace{-1em}%
    \caption{Ignoring memory issues, off-the-shelf methods require intractable computation time to compute a single
high-dimensional matrix exponential for the 3D Heat Diffusion system used in our evaluation.}%
\label{fig:mat_exp}%
\vspace{-1.5em}
\end{figure}

As time advances, the behaviors of affine systems can be written in terms of the matrix exponential,
which can be used to compute the reachable set.
For high-dimensional systems, however, computing the matrix exponential becomes both the runtime and memory bottleneck.
A simple experiment demonstrating the runtime problem is shown in Figure~\ref{fig:mat_exp}, where extrapolation
predicts it would take over a year to compute a single matrix exponential once a system has over ten million dimensions.
A second problem is memory: although the $A$ matrix for large systems can often be sparse,
the matrix exponential of $A$ is dense.
%
%
The amount of memory needed to simply store the matrix exponential result, a dense $n \times n$ matrix,
can greatly exceed what is available, even if its computation time was instant.
For a million-dimensional system, this matrix would have $10^{12}$ numbers and need about $8$ TB of main memory.


The main contribution of this paper is an approach for affine systems verification that can,
with sufficient problem structure, scale to
extremely high-dimensional systems, thousands to millions of times larger than existing methods.
We overcome the memory and computation time problems through an original \emph{combination} of new and existing techniques.
The memory improvements are possible through a method that uses both aspects of reachability with
support functions~\cite{le2010reachability} (projecting the reachable set onto an output space) and affine representations~\cite{han2006reachability}
/ zonotopes~\cite{girard2005reachability} (low-dimensional initial spaces).
The run-time improvements use simulations to compute parts of the
matrix exponential~\cite{duggirala2016parsimonious}.
%
Since large dynamics matrices are often sparse (and must be sparse to simply fit into memory), we can perform numerical
simulations using efficient Krylov subspace methods~\cite{gallopoulos1992efficient,han2006reachability}.
We use a recently-published a posteriori error bound~\cite{wang2017error} to determine when the dimension
of the Krylov subspace is sufficient for an accurate simulation result.
Further, we propose a modification to the Lanczos iteration (used by Krylov subspace methods) that reduces memory requirements, allowing us to
significantly increase the number of iterations before memory is exhausted.
Although some of the techniques have been used individually before, no existing method for affine systems
has demonstrated scalability beyond a few thousand dimensions.
In our evaluation, the largest system we analyze has one billion dimensions.

\new{
  The research presented here builds off a previous workshop paper~\cite{bak18arch}, with several new developments:
  (i) this paper focuses on scalability in high dimensions (up to $10^9$ dimensions), whereas the earlier work performed a runtime improvement evaluation on comparatively smaller systems (up to $10^4$ dimensions);
  (ii) Section~\ref{sec:computational} provides a detailed description and pseudo-code for the modified Krylov methods, which use an a posteriori error bound (Lemma~\ref{lm:error-bound}), as opposed to the relative error estimate in the earlier paper;
  (iii) we provide memory improvements to the projected Lanczos iteration in Section~\ref{ssec:lanczos} which, in our evaluation in Section~\ref{ssec:heat3d}, is shown to increase scalability by two orders of magnitude.
}

Section~\ref{sec:basic} first reviews affine discrete-time safety verification,
which uses an $n \times n$ matrix exponential at each time step in the analysis.
Next, Section~\ref{sec:memory} presents memory improvements followed by Section~\ref{sec:computational},
which focuses on reducing computation time.
%
%
An evaluation on several large benchmarks, including a 3D Heat Diffusion system with one billion dimensions,
is given in Section~\ref{sec:evaluation}, followed by a review of related work and a conclusion.

%% file: basic.tex
\section{Affine Verification Review}
\label{sec:basic}

An affine, discrete-time, bounded safety verification problem \old{is provided as input} \new{is defined by} the system dynamics
$\dot{x} = Ax + b$, a set of initial states $\mathcal{I}$ defined as all states $x_0$ where
the linear constraints $\mathcal{I}_x x_0 \leq \iota_x$ hold, unsafe states $\mathcal{U}$ defined with linear
constraints $\mathcal{U}_x x \leq \upsilon_x$, a step size $\delta$ and time bound $T$.
The system is called unsafe if and only if there exists a time $t = k \delta \leq T$ such that
$x_0 \in \mathcal{I}$, $x = e^{At} x_0$, and $x \in \mathcal{U}$.
The goal is to prove a system is safe or find a counter-example, which can be defined by an initial state $x_0$
and time $t$.

\subsection{Basic Verification Approach}
\label{ssec:basic_approach}

An affine system with dynamics $\dot{x} = Ax + b$ can be verified by first converting it to a
linear system (without the $b$ term), by adding a fresh variable to account for the effects of the forcing term $b$.
The new $A$ matrix has an extra column consisting of the entries of the $b$ vector, and an extra row of all zeros.
The initial value of the new variable is assigned to $1$, and, since the row in $A$ defining its differential equation
is all zeros, the new variable's value remains at $1$ at all times.
Thus, the effect of the extra column in the $A$ matrix is the same as the $b$ vector
in the original system.
We consider linear systems after
this transformation, assuming the form $\dot{x} = Ax$.

Safety can be checked by constructing a LP at each discrete time $t$ that contains two copies of the state variables,
$x_0$ and $x$, and encodes the initial state conditions $\mathcal{I}_x x_0 \leq \iota_x$, the unsafe state
conditions $\mathcal{U}_x x \leq \upsilon_x$, and the linear relationship (for a fixed $t$) between the initial and
final variables $x = e^{At} x_0$.
If the LP is feasible, the solver provides an assignment to the variables that can be used to
construct the counter-example.
The bulk of the computation time is spent on these two operations: (i) computing $e^{At}$ and (ii) solving the LP.

\subsection{Timed Harmonic Oscillator Example}
\label{ssec:ha}

We will use \old{am} \new{an} example of a timed harmonic oscillator to demonstrate the methods in this paper.
The timed harmonic oscillator is a system with dynamics $\dot{x} = y$, $\dot{y} = -x$, and $\dot{t} = 1$.
For the initial set of states, take $x_0 = -5$, $y_0 \in [0, 1]$, and $t_0 = 0$.
The unsafe set of states consists of all states where $x = 4$.
We attempt to verify the system with a discrete time step of $\delta=\frac{\pi}{4}$ and a time bound of $T=\pi$.

\begin{figure}[t]%
    \centering%
    \includegraphics[width=0.7\textwidth]{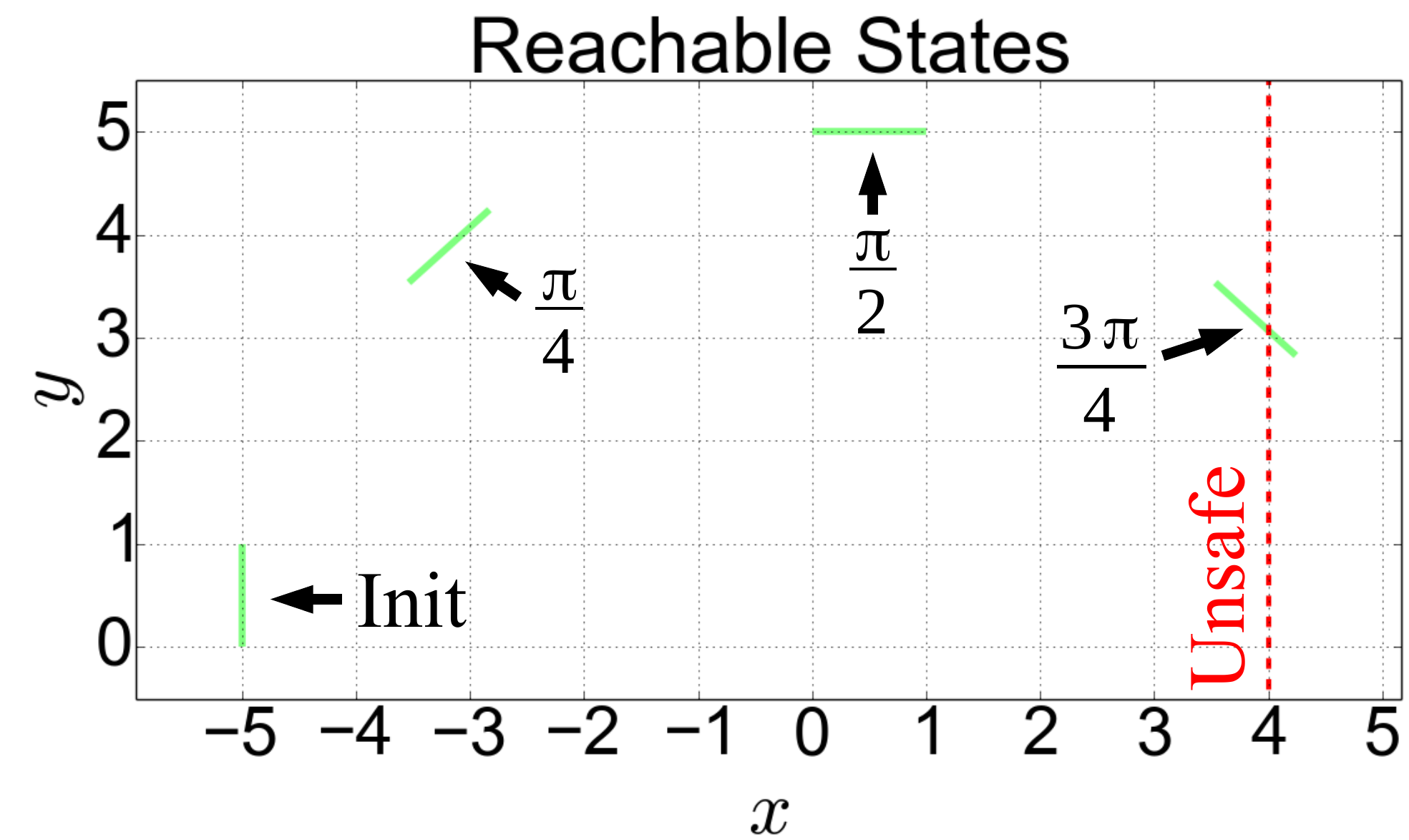}%
		\vspace{-1em}%
    \caption{The timed harmonic oscillator system can reach an unsafe state at
time $\frac{3\pi}{4}$.}%
\vspace{-1em}%
\label{fig:ha_plot}%
\vspace{-0.5em}%
\end{figure}

On the $x$-$y$ plane, solutions of the system rotate clockwise around the origin.
The reachable set is shown in Figure~\ref{fig:ha_plot}.
From the figure, it is apparent that at time $\frac{3 \pi}{4}$, the unsafe states are reachable.

We can show this computationally and find the initial state that leads to the violation.
First, we convert the three-variable affine system (the $t$ variable has an affine term),
to a four-variable linear system
using the affine-to-linear transformation described before.
The variables in the transformed system are $\vec{x} = (x, y, t, a)^T$, where $a$ is the
newly-introduced variable, which is initially 1 and remains constant at all times.
The transformed dynamics now form the four-dimensional linear system $\dot{\vec{x}} = A \vec{x}$, with
\begin{equation*}
A =
\begin{pmatrix}
   0 & 1 & 0 & 0 \\
   -1 & 0 & 0 & 0 \\
   0 & 0 & 0 & 1 \\
   0 & 0 & 0 & 0
   \end{pmatrix}.
\end{equation*}

\input{outline_fig}

Next, at each discrete time step, we construct a set of linear constraints.
The constraints have two copies of the variables, one set at the initial time
$\{x_0, y_0, t_0, a_0\}$, and one set at the current time step $\{x, y, t, a\}$.
The linear constraints at the time of the violation, $\frac{3\pi}{4}$, are shown
in Figure~\ref{sfig:full_constraints}.
The only constraints that change between time steps are the ones encoding the matrix exponential at the
current time (the values surrounded by a red rounded rectangle), for which we
reuse terminology~\cite{bak2017arch} and refer to as the \emph{basis matrix}.

The linear constraints are then passed to an LP solver to check if they are feasible.
For time steps $0$, $\frac{\pi}{4}$, and $\frac{\pi}{2}$, the LP solver returns that no solution exists.
At time $\frac{3\pi}{4}$, which has the constraints shown in the figure, the LP solver finds a feasible solution,
and provides an assignment to the variables.
In particular, its output indicates that starting from initial state
$(x_0, y_0, t_0, a_0)^T$ = $(-5, 0.66, 0, 1)^T$,
the system can reach the unsafe state $(x, y, t, a)^T$ = $(4, 3.07, 2.36, 1)^T$.

%% file: outline_fig.tex
\definecolor{g}{rgb}{0.0,0.7,0.0}

\begin{figure*}[t] %
\centering
\subfloat[Basic approach with full linear constraints (Section~\ref{ssec:basic_approach})]{
\label{sfig:full_constraints}
\resizebox{0.75\columnwidth}{!}{
\input{outline1}
}
} \\
\subfloat[Projecting onto the output space (Section~\ref{sec:projecting_key})]{
\label{sfig:projecting_key}
\resizebox{0.5\columnwidth}{!}{
\input{outline2}
}
} ~ ~~
\subfloat[Projecting from the initial space onto the output space (Section~\ref{sec:combining_fixed})]{
\label{sfig:combining_fixed}
\resizebox{0.48\columnwidth}{!}{
\input{outline3}
}
}%
\caption{The linear constraints at time $\frac{3\pi}{4}$
for the timed harmonic oscillator example described in Section~\ref{ssec:ha} can be encoded in
different ways. At each time step, only the basis matrix changes in the constraints.}%
\label{fig:ha_constraints}%
\end{figure*}

%% file: outline1.tex
\begin{tikzpicture}[baseline=(m.west)]
\matrix (m) [
         left delimiter={(},right delimiter={)}, matrix of math nodes, ampersand replacement=\&, column sep={3em,between origins}]
         {
-1 \& 0  \& 0  \& 0  \& -0.707~~~~ \& 0.707 \& 0  \& 0  \\ 
0  \& -1 \& 0  \& 0  \& -0.707~~~~ \& -0.707 \& 0  \& 0  \\ 
0  \& 0  \& -1 \& 0  \& 0     \& 0     \& 1  \& 2.36 \\
0  \& 0  \& 0  \& -1 \& 0     \& 0     \& 0  \& 1 \\
|[blue, fill=blue!20]| \bf{1}  \& |[blue, fill=blue!20]| \bf{0}  \& |[blue, fill=blue!20]| \bf{0}
\& |[blue, fill=blue!20]| \bf{0} \& 0 \& 0 \& 0  \& 0 \\
0  \& 0  \& 0  \& 0  \& |[g]| \bf{1}  \& |[g]| \bf{0}  \& |[g]| \bf{0}  \& |[g]| \bf{0} \\
0  \& 0  \& 0  \& 0  \& |[g]| \bf{0}  \& |[g]| \bf{-1} \& |[g]| \bf{0}  \& |[g]| \bf{0} \\
0  \& 0  \& 0  \& 0  \& |[g]| \bf{0}  \& |[g]| \bf{1}  \& |[g]| \bf{0}  \& |[g]| \bf{0} \\
0  \& 0  \& 0  \& 0  \& |[g]| \bf{0}  \& |[g]| \bf{0}  \& |[g]| \bf{1}  \& |[g]| \bf{0} \\
0  \& 0  \& 0  \& 0  \& |[g]| \bf{0}  \& |[g]| \bf{0}  \& |[g]| \bf{0}  \& |[g]| \bf{1} \\
};
\draw [red,ultra thick,rounded corners]
  ($(m.south west) !.81! (0.6, 1.0)$)
  rectangle
  ($(m.south east) !0.98! (m.north east)$);
\draw [] (1.8, 2.4) node [red]
{\bf{Basis Matrix}};
\draw [align=center] (-2.0, -0.5) node [blue, fill=blue!20]
{\bf{Unsafe} \\ \bf{Condition}};
\draw [g,double]
  ($(m.south west) !.81! (0.6, -2.05)$)
  rectangle
  ($(m.south east) !0.49! (m.north east)$);
\draw [] (1.8, -2.3) node [g]
{\bf{Initial State Conditions}};
\end{tikzpicture}
\begin{tikzpicture}[baseline=(m.west)]
\matrix (m) [left delimiter={(},right delimiter={)}, matrix of math nodes, row sep={1.6em,between origins}]
{
x \\
y \\
t \\
a \\
x_0 \\
y_0 \\
t_0 \\
a_0 \\
};
\end{tikzpicture}
\begin{tikzpicture}[baseline=(m.west), column 1/.style={anchor=base west}]
\matrix (m) [matrix of math nodes]
{=~~0\\
=~~0\\
=~~0\\
=~~0\\
|[blue, fill=blue!20]| \bf{=~~4} \\
|[g]| \bf{=-5}\\
|[g]| \bf{\leq~~ 0}\\
|[g]| \bf{\leq~~1}\\
|[g]| \bf{=~~0}\\
|[g]| \bf{=~~1}\\
};
\draw [g,double]
  ($(m.south west) !0.0! (m.north west)$)
  rectangle
  ($(m.south east) !0.5! (m.north east)$);
\end{tikzpicture}

%% file: outline2.tex
\begin{tikzpicture}[baseline=(m.west)]
\matrix (m) [
         left delimiter={(},right delimiter={)}, matrix of math nodes, ampersand replacement=\&, column sep={3em,between origins}]
         {
-1~~  \& -0.707~ \& ~0.707 \& 0  \& 0  \\ 
|[blue, fill=blue!20]| \bf{1}  \& 0 \& 0 \& 0  \& 0 \\
0  \& |[g]| \bf{1}  \& |[g]| \bf{0}  \& |[g]| \bf{0}  \& |[g]| \bf{0} \\
0  \& |[g]| \bf{0}  \& |[g]| \bf{-1} \& |[g]| \bf{0}  \& |[g]| \bf{0} \\
0  \& |[g]| \bf{0}  \& |[g]| \bf{1}  \& |[g]| \bf{0}  \& |[g]| \bf{0} \\
0  \& |[g]| \bf{0}  \& |[g]| \bf{0}  \& |[g]| \bf{1}  \& |[g]| \bf{0} \\
0  \& |[g]| \bf{0}  \& |[g]| \bf{0}  \& |[g]| \bf{0}  \& |[g]| \bf{1} \\
};
\draw [] (0.45, 1.73) node [red]
{\bf{Basis Matrix}};
\draw [red,ultra thick,rounded corners]
  ($(m.south west) !.83! (-1.3, 1.6)$)
  rectangle
  ($(m.south east) !0.98! (m.north east)$);
\draw [g,double]
  ($(m.south west) !0.8! (-1.3, 1.1)$)
  rectangle
  ($(m.south east) !0.05! (m.north east)$);
\draw [] (0.45, -1.6) node [g]
{\bf{Initial State Conditions}};
\end{tikzpicture}
\begin{tikzpicture}[baseline=(m.west)]
\matrix (m) [left delimiter={(},right delimiter={)}, matrix of math nodes, row sep={1.6em,between origins}]
{
o_x \\
x_0 \\
y_0 \\
t_0 \\
a_0 \\
};
\end{tikzpicture}
\begin{tikzpicture}[baseline=(m.west), column 1/.style={anchor=base west}]
\matrix (m) [matrix of math nodes]
{=~~0\\
|[blue, fill=blue!20]| \bf{=~~4} \\
|[g]| \bf{=-5}\\
|[g]| \bf{\leq~~ 0}\\
|[g]| \bf{\leq~~1}\\
|[g]| \bf{=~~0}\\
|[g]| \bf{=~~1}\\
};
\draw [g,double]
  ($(m.south west) !0.03! (m.north west)$)
  rectangle
  ($(m.south east) !0.68! (m.north east)$);
\end{tikzpicture}

%% file: outline3.tex
\begin{tikzpicture}[baseline=(m.west)]
\matrix (m) [
         left delimiter={(},right delimiter={)}, matrix of math nodes, ampersand replacement=\&, column sep={3em,between origins}]
         {
-1  \& 0.707~~~~ \& 3.54 \\ 
|[blue, fill=blue!20]| \bf{1}  \& 0 \& 0 \\
0  \& |[g]| \bf{-1} \& |[g]| \bf{0} \\
0  \& |[g]| \bf{1}  \& |[g]| \bf{0} \\
0  \& |[g]| \bf{0}  \& |[g]| \bf{1} \\
};
\draw [] (0.4, 1.3) node [red]
{\bf{Basis Matrix}};
\draw [red,ultra thick,rounded corners]
  ($(m.south west) !.81! (-0.5, 1.05)$)
  rectangle
  ($(m.south east) !0.97! (m.north east)$);
\draw [] (0.4, -1.25) node [g]
{\bf{Initial State Conditions}};
\draw [g, double]
  ($(m.south west) !.81! (-0.5, 0.55)$)
  rectangle
  ($(m.south east) !0.05! (m.north east)$);
\end{tikzpicture}
\begin{tikzpicture}[baseline=(m.west)]
\matrix (m) [left delimiter={(},right delimiter={)}, matrix of math nodes, row sep={1.6em,between origins}]
{
o_x \\
i_y \\
i_f \\
};
\end{tikzpicture}
\begin{tikzpicture}[baseline=(m.west), column 1/.style={anchor=base west}]
\matrix (m) [matrix of math nodes]
{=~~0\\
|[blue, fill=blue!20]| \bf{=~~4} \\
|[g]| \bf{\leq~~ 0}\\
|[g]| \bf{\leq~~1}\\
|[g]| \bf{=~~1}\\
};
\draw [g,double]
  ($(m.south west) !0.03! (m.north west)$)
  rectangle
  ($(m.south east) !0.575! (m.north east)$);
\end{tikzpicture}

%% file: memory.tex
\section{Memory Improvements}
\label{sec:memory}

Although the basic verification approach works, it does not scale to very high dimensions.
As mentioned in the introduction, computing and storing the basis matrix (the matrix exponential)
is typically the bottleneck to verification scalability.
In this section we focus on the memory problem, and show how we can reduce the height (Section~\ref{sec:projecting_key}) 
and width (Section~\ref{sec:combining_fixed}) of the basis matrix, by taking advantage of common problem structure.
%

\subsection{Projecting onto the Output Space}
\label{sec:projecting_key}

First, we reduce the \emph{height} of the basis matrix 
(compare the basis matrices in Figure~\ref{sfig:full_constraints} and Figure~\ref{sfig:projecting_key}).
This is done by a method similar to the use of support functions with a fixed number of directions for 
reachability analysis~\cite{le2010reachability}.
The common problem structure exploited is that the verification result often only depends on a small number of 
directions, much smaller than the number of system variables.

Depending on the type of problem being solved (linear verification, plotting, or hybrid automaton reachability), 
these directions arise from different sources.
For a safety verification problem for linear systems, these directions come from each of the constraints in the 
conjunction defining the unsafe states.
For a plot, we only need to compute a projection onto the two or three plot dimensions.
In this case, the important directions are the unit vectors in each of these dimensions.
Plots can then be produced efficiently by running multiple optimizations over projections 
of the convex reachable set at each time step~\cite{ray16vertex_enumeration,lotov2013interactive}.
For the hybrid automaton setting, additional directions can come from the constraints in the mode 
invariants, as well as from the guard conditions.

We can combine these directions into an output matrix $C$, where the output variables are $y = Cx$, and the height of the 
matrix is the number of output directions, $o$, needed for the current problem.
The unsafe states, $\mathcal{U}$, are then redefined in the output space, $\mathcal{U}_y y \leq \upsilon_y$.
Finally, the basis matrix in the constraints is the $o \times n$ projection of the matrix exponential 
onto the output space, $C e^{At}$.

Consider applying this approach to the timed harmonic oscillator system 
of Section~\ref{ssec:ha}, where the unsafe states are defined by $x = 4$.
The other three dimensions, $y$, $t$, and $a$, do not impact the result of the safety check, and so they 
(and their corresponding constraints) can be removed from the set of linear constraints, 
as is done in Figure~\ref{sfig:projecting_key}.
In this case, the output matrix for this system is the $1 \times 4$ matrix 
$C = \begin{psmallmatrix} 1 & 0 & 0 & 0 \end{psmallmatrix}$.
We then define the unsafe states in terms of the single output space variable $o_x$, and replace the 
basis matrix by the projected matrix exponential $C e^{At}$.


\subsection{Projecting from the Initial Space}
\label{sec:combining_fixed}

Next, we reduce the \emph{width} of the basis matrix (compare the basis matrix in Figure~\ref{sfig:projecting_key} 
and Figure~\ref{sfig:combining_fixed}).
This is done with a method similar to reachability using affine representations~\cite{han2006reachability}, and
is also similar to reachability with zonotopes~\cite{girard2005reachability} with a small number of generators.
The common problem structure exploited is that the initial states are often low-dimensional.
For example, there may not be uncertainty in every variable, or the initial states of variables may be related.

In this case, we can define an $i$-dimensional initial space using an $n \times i$ matrix $E$, where
the initial states $z$ are related to the original variables by $x = Ez$.
The initial states $\mathcal{I}$ are then redefined with constraints in the initial space,
$\mathcal{I}_z z \leq \iota_z$.
The $o \times i$ basis matrix is now computed using both projections, $C e^{At} E$.

In the timed harmonic oscillator system of Section~\ref{ssec:ha}, we can define the initial states 
using $i=2$ dimensions. 
These are $i_y$, which corresponds to the initial $y$ value, 
and $i_f$ which is the fixed initial values of all the other variables.
The $E$ matrix is the $4 \times 2$ matrix $\begin{psmallmatrix} 0 & 1 & 0 & 0 \\ -5 & 0 & 0 & 1 \end{psmallmatrix}^T$,
and the initial constraints are $0 \leq i_y \leq 1$ and $i_f = 1$.
The basis matrix is the product $C e^{At} E$ at each step.

Using both methods, we have reduced the basis matrix from an $n \times n$ matrix to a $o \times i$ matrix.
Importantly, we do not need both $o$ and $i$ to be very small for this reduction to be useful, only their \emph{product}.
Given, say, 800 MB to store the basis matrix ($10^8$ double-precision numbers), 
the original approach would fill the memory when $n = 10^4$, a ten-thousand dimensional system.
In contrast, a million-dimensional system with every dimension initially independent and uncertain, $i=10^6$, could 
still be analyzed as long as the unsafe states are defined using less than $100$ output directions.

%% file: computational.tex
\section{Computation Time Improvements}
\label{sec:computational}

Although we can define the smaller basis matrix using $C e^{At} E$, this does not help in terms of computation time if
we use the direct approach of computing $e^{At}$ at each step and then multiplying by $C$ and $E$.
In this section, we describe a series of improvements targeting the computational efficiency of the method.

\subsection{Basis Matrix using Numerical Simulations}
\label{ssec:basis_with_simulation}

There are many ways to compute the matrix exponential~\cite{moler2003nineteen}.
Generally, the methods implemented in off-the-shelf libraries use a combination of squaring and scaling and Pade
approximation~(methods $2$ and $3$~\cite{moler2003nineteen}), which compute the entire matrix at once.

Here, we instead use an alternative method to compute the matrix exponential, using a series of numerical simulations
(method $5$~\cite{moler2003nineteen}).
The matrix exponential is computed one column at time by using the fact that
$e^{At} = e^{At} ~ \mathbf{I}_{n \times n} = e^{At} ~ (\mathbf{e}_1 | \mathbf{e}_2 | \ldots | \mathbf{e}_n)$.
The $j$th column of $e^{At}$ is equal to $e^{At} \mathbf{e}_j$, where $\mathbf{e}_j$ is
the $j$th column of the identity matrix.
The value of $e^{At} \mathbf{e}_j$, however, is just the solution of the linear system $\dot{x} = Ax$ at time $t$ from
initial state $x(0) = \mathbf{e}_j$.
To compute this, we can perform a numerical simulation with an off-the-shelf numerical method such as Runge-Kutta.
This process is repeated for each column of the identity matrix to compute the full matrix
exponential.
For the verification problem, we need the value of $e^{At}$ at multiple time steps, and so we run
the numerical simulations up to the time bound $T$, recording the value at each multiple of the step size $\delta$.
The values from each column are then combined at each multiple of time step to form the
basis matrix in the LP~\cite{duggirala2016parsimonious}.

We propose to adapt this method to take advantage of initial and output spaces.
Since we need to compute the basis matrix $C e^{At} E$, rather than using each column of the identity matrix,
we can instead compute simulations from each column of the $E$ matrix, and then project the state in the simulation using
the $C$ matrix.
We compute $e^{At} E$ by noting that
$e^{At} ~ E = e^{At} ~ (E_{*,1} | E_{*,2} | \ldots | E_{*,i})$,
where $E_{*,j}$ is the $j$th column of $E$.
As before with the identity matrix, each column $j$ can be computed separately with a numerical simulation
of the linear system $\dot{x} = Ax$ at time $t$ from initial state $x(0) = E_{*,j}$.
There are $i$ columns in $E$, corresponding to the $i$ dimensions of the initial states.
If $i$ is much smaller than $n$, this approach will be significantly faster than
computing the full matrix exponential and then doing the multiplication with $E$.
For the timed-harmonic oscillator system constraints in Figure~\ref{sfig:combining_fixed}, for example, since the
dimension of the initial space $i=2$, the basis matrix could be computed in this fashion using two numerical
simulations.

If the initial state dimension is large, the computation may still require a large number of simulations.
In this work, we propose \new{a} new approach that can reduce the required number of simulations if the
output space is small.
The method works by performing simulations using the transpose system dynamics.
Since $C e^{At} E = ((C e^{At} E)^T)^T = (E^T (e^{At})^T C^T)^T = (E^T e^{A^Tt} C^T)^T$, the basis matrix can also
be computed by performing $o$ simulations (one for each column of $C^T$).
We compute $e^{A^Tt} ~ C^T$ by noting that
$e^{A^Tt} ~ C^T = e^{A^Tt} ~ (C_{*,1}^T | C_{*,2}^T | \ldots | C_{*,o}^T)$,
and performing a numerical simulation
of the linear system $\dot{x} = A^T x$ up to time $t$ from initial state $x(0) = C_{*,j}^T$, for each column $j$ of $C^T$.
The results are then multiplied by $E^T$, and transposed to recover the basis matrix.
This allows us to compute values of the basis matrix one \emph{row} at a time, and so we can compute the basis matrix
using only $o$ numerical simulations.
In practice, only one of these is necessary, and so we can choose the minimum of $i$ and $o$ and perform that
many numerical simulations, rather than computing an $n \times n$ matrix exponential.

In the timed-harmonic oscillator system, for example, since $o=1$, the entire basis matrix at each step can be
computed with a single numerical simulation.
Starting from the state corresponding to the single
output direction $(1, 0, 0, 0)^T$, we can simulate
using the transpose dynamics $A^T$ up to time $\frac{3\pi}{4}$ to get the state $(-0.707, 0.707, 0, 0)$.
This is then
projecting with $E^T = \begin{psmallmatrix} 0 & 1 & 0 & 0 \\ -5 & 0 & 0 & 1 \end{psmallmatrix}$ to get
$(0.707, 3.54)^T$, which is transposed to get the basis matrix in
Figure~\ref{sfig:combining_fixed}.

\subsection{Simulations using the Krylov Subspace}
\label{ssec:krylov}

When the system matrix $A$ is high-dimensional, it is also often sparse (in fact, if $A$ has more than tens of
thousands of dimensions and can fit in memory, it must be sparse or otherwise compressed).
We can exploit this structure to speed up numerical simulations.

The Krylov subspace simulation method~\cite{gallopoulos1992efficient} computes an approximation
of $e^{A}v$, where $v$ is some initial state.
This is done by finding the element of the $k$-dimensional Krylov subspace
$K_{k} \equiv span\{v,Av,\dots,A^{k-1}v\}$ that best approximates $e^{A}v$.
\new{Intuitively, the $k$-dimensional Krylov subspace can exactly represent the first $k$ terms of the Taylor expansion of}
$e^{A}v$, \new{making it a good candidate for accurate approximation.}
We do not review the full theory here, but instead focus on computational aspects as they
relate to the verification problem.

Note, however, that each simulation has a different initial state, and different initial states $v$ will have different Krylov subspaces.
This is important, since it means we are \textbf{not} proposing to verify the system through an abstraction of the dynamics
matrix $A$ by a single lower-dimensional system matrix.

The approximation uses a fixed number of iterations of the
well-known Arnoldi algorithm~\cite{arnoldi1951principle,trefethen1997numerical}.
The pseudocode is shown in Algorithm~\ref{alg:arnoldi}.
The Arnoldi algorithm computes an orthonormal basis for the Krylov subspace $K_{k}$ by
starting with a normalized version of $v$ as both the first orthonormal direction and the current vector and,
at each iteration, (1) multiplying the current vector by $A$ (line~\ref{line:arnoldi_mult}),
(2) projecting out the previous orthonormal directions from the current vector
(the loop on lines~\ref{line:arnoldi_proj_start}-\ref{line:arnoldi_proj_end}),
(3) normalizing the current vector (lines~\ref{line:arnoldi_normalize}-\ref{line:arnoldi_add}),
and (4) adding it to the list of orthonormal directions (line~\ref{line:arnoldi_add}).
If the norm computed on line~\ref{line:arnoldi_normalize} is ever zero,
the loop can terminate early (not shown) and the approximation will be exact.
The memory needed for the Arnoldi iteration, which can be obtained by looking at the sizes of the outputs, is thus:
\begin{equation}
\label{eq:memory_arnoldi}
k \times (n + k) \times \texttt{sizeof(double)}
\end{equation}

\algnewcommand{\algorithmicredif}{\red{\textbf{if}}}
\algnewcommand{\algorithmicredthen}{\red{\textbf{then}}}
\algnewcommand{\algorithmicredend}{\red{\textbf{end}}}
\algnewcommand{\algorithmicredelse}{\red{\textbf{else}}}

\algdef{SE}[REDIF]{RedIf}{RedEndIf}[1]{\algorithmicredif\ #1\ \algorithmicredthen}{\algorithmicredend\ \algorithmicredif}%
\algdef{Ce}[REDELSE]{REDIF}{RedElse}{RedEndIf}{\algorithmicredelse}
\algtext*{RedEndIf} 

\begin{algorithm}[t]
  \caption{Original Arnoldi algorithm}\label{alg:arnoldi}
  \begin{algorithmic}[1]
\Require{normalized init $n \times 1$ vector $v$, $n \times n$ matrix $A$, iterations $k$}
\Ensure{$n \times k$ matrix $V$, $k \times k$ matrix $H$}
\State $V_{*,1} \gets v$ \Comment{assign to first column of $V$}
\For{$i$ \textbf{from} $2$ \textbf{to} $k+1$}
\State $V_{*,i} \gets A ~ V_{*,i-1}$ \label{line:arnoldi_mult}
\For{$j$ \textbf{from} $1$ \textbf{to} $i$} \label{line:arnoldi_proj_start}
\State $H_{j,i-1} \gets (V_{*,j})^T V_{*,i}$
\State $V_{*,i} \gets V_{*,i} - H_{j,i-1} V_{*,j}$ \label{line:arnoldi_proj_end}
\EndFor
\State $H_{i,i-1} \gets ||V_{*,i}||$ \label{line:arnoldi_normalize}
\State $V_{*,i} \gets \frac{V_{*,i}}{H_{i,i-1}}$ \label{line:arnoldi_add}
\EndFor
\State $H \gets H_{1:k, *}$ \Comment{discard extra row of $H$}
\State $V \gets V_{*, 1:k}$ \Comment{discard extra column of $V$}
  \end{algorithmic}
\end{algorithm}

After $k$ iterations complete, the outputs are two matrices $V$ and $H$, which we refer to as $V_k$ and $H_k$.
$V_k$ is the $n \times k$ matrix of orthonormal basis vectors
and $H_k$ is the $k \times k$ matrix that is a projection of the linear transformation $A$ in the Krylov subspace $K_{k}$.

The outputs of the Arnoldi algorithm can be used to approximate $e^{A}v$.
This is done by projecting the initial $n$-dimensional state onto the smaller, $k$-dimensional Krylov subspace, computing the matrix exponential
using the projected linear transformation $H_k$,
and then projecting the result back to the original $n$-dimensional space using $V_k$.
By the design of the Krylov subspace, the projection of the initial state $v$ is just
the length of $v$ multiplied by the first unit vector in the subspace, $\mathbf{e}_1$.
Further, since for any time $t$, the Krylov subspaces associated with
$A$ and $At$ are identical \new{(because $span\{v,Av,\dots,A^{k-1}v\}$ is the same as $span\{tv,tAv,\dots,tA^{k-1}v\}$)},
we can use the same $V_k$ and $H_k$ to compute the approximation at any point in time.
The formula for the approximation is:
\begin{equation}
\label{eq:approx-formula}
e^{At}v \approx \norm{v} V_{k} e^{H_{k}t} \mathbf{e}_1
\end{equation}

Equation~\ref{eq:approx-formula} is especially useful when the size of $A$ is huge, e.g., millions of dimensions,
since it transforms the computation with a large matrix $A$ to a problem
with a much smaller matrix $H_{k}$.
For fast computation, we would like to minimize the size of $H_{k}$ by using a small number
of Arnoldi iterations $k$, but this has the effect of reducing the approximation accuracy.
Thus, it is critical to select $k$ large enough to be accurate, but small enough to be fast.

Earlier work on reachability with Krylov subspace methods~\cite{han2006reachability} used an a priori error
bound~\cite{gallopoulos1992efficient} to determine $k$.
The error of the approximation for a fixed $k$ is bounded by
\begin{equation}
\label{eq:error-bound}
 \norm{v} \frac{\norm{At}^{k}e^{\norm{At}}}{k!}.
\end{equation}

Unfortunately, the a priori error bound can often be unusably pessimistic.
For example, one of the models we will use in our evaluation is a 100x100x100 3D Heat Diffusion system (one million dimensions).
At time 50, this system has matrix norm $\norm{At} = 32771611$.
For an initial unit vector with $\norm{v} = 1$, even using a full dimensional Krylov subspace ($k=10^6$),
the computed a priori error bound from Equation~\ref{eq:error-bound} is effectively unusable, $10^{16182319}$.

\begin{algorithm}[t]
  \caption{Arnoldi algorithm with a posteriori error control}\label{alg:arnoldi_error}
  \begin{algorithmic}[1]
\Require{normalized init $n \times 1$ vector $v$, $n \times n$ matrix $A$, \red{error target $\epsilon$}}
\Ensure{$n \times k$ matrix $V$, $k \times k$ matrix $H$}
\State $V_{*,1} \gets v$
\State \red{$k \gets 4$}
\For{$i$ \textbf{from} $2$ \textbf{to} \red{$\infty$}}
\State $V_{*,i} \gets A ~ V_{*,i-1}$
\For{$j$ \textbf{from} $1$ \textbf{to} $i$}
\State $H_{j,i-1} \gets (V_{*,j})^T V_{*,i}$
\State $V_{*,i} \gets V_{*,i} - H_{j,i-1} V_{*,j}$
\EndFor
\State $H_{i,i-1} \gets ||V_{*,i}||$
\State $V_{*,i} \gets \frac{V_{*,i}}{H_{i,i-1}}$
\RedIf{\red{$i = k$}} \red{\Comment{check error upon reaching $k$ iterations}}
\RedIf{\red{$\texttt{compute-error}(A, H_{1:i-1,*}) < \epsilon$}}
    \State \red{\textbf{break}}
    \RedElse
    \State \red{$k \gets \texttt{ceil}(1.1 * k)$}
\RedEndIf
\RedEndIf
  \EndFor
\State $H \gets H_{1:k, *}$ \Comment{discard extra row of $H$}
\State $V \gets V_{*, 1:k}$ \Comment{discard extra column of $V$}
  \end{algorithmic}
\end{algorithm}

In this work, we instead use a recently-developed a posteriori error bound~\cite{wang2017error}, which uses information
from the $H$ matrix as well as the extreme eigenvalues of $A$ to compute a bound on the error.
\new{The bound works by creating an error function using the log norm of $A$ and looking at the derivative of this error over time.}
\begin{lemma}[\cite{wang2017error}]\label{lm:error-bound}
Let $A \in \mathbb{R}^{n \times n}$ and $v \in \mathbb{R}^n$ with $\norm{v} = 1$. Let $V_k$ be the orthogonal matrix and $H_k$ be the upper Hessenberg matrix generated by the Arnoldi process for A and $v$. Let $\omega_k(\tau) = V_ke^{-\tau H_k}e_1$ be the Arnoldi approximation to $\omega(\tau) = e^{-\tau A}v$. Then the approximation error satisfies
\begin{equation}
\norm{\omega(\tau) - \omega_k(\tau)} \leq h_{k+1, k}e^{-min\{\nu(A), 0\}\tau}\int_{0}^{\tau}|h(t)|dt,
\end{equation}
where $h(t) := e^T_ke^{-tH_k}e_1$ is the $(k, 1)$ entry of the matrix $e^{-tH_k}$ and $\nu(A)$ is the smallest eigenvalues of $(\frac{A + A^T}{2})$.
\end{lemma}
The above lemma computes the error bound of approximating $\omega(\tau) = e^{-\tau A}v$ with the Arnoldi algorithm.
In our application, we want to approximate $e^{\tau A}v = e^{-\tau (-A)}v$.
To do that, we only need to feed $B = -A$ as an input to the Arnoldi algorithm
and use the lemma with the matrix $B$ when computing the error bound.

Since this error bound uses values of the $H$ matrix which is an output of the Arnoldi algorithm,
we cannot determine $k$ ahead of time, as we could with an a priori bound.
However, with this bound we can provide an accuracy guarantee with significantly fewer iterations.

\begin{algorithm}[t]
  \caption{Original Lanczos algorithm}\label{alg:lanczos}
  \begin{algorithmic}[1]
\Require{normalized init $n \times 1$ vector $v$, $n \times n$ matrix $A$, iterations $k$}
\Ensure{$n \times k$ matrix $V$, $k \times k$ matrix $H$}
\State $V_{*,1} \gets v$
\For{$i$ \textbf{from} $2$ \textbf{to} $k+1$}
\State $V_{*,i} \gets A ~ V_{*,i-1}$

\If{i > 2} \label{line:lanczos_j=i-2_start}
\State $H_{i-2,i-1} \gets H_{i-1,i-2}$
\State $V_{*,i} \gets V_{*,i} - H_{i-2,i-1} V_{*,i-2}$  \label{line:lanczos_j=i-2_end}
\EndIf
\State $H_{i-1,i-1} \gets (V_{*,i-1})^T V_{*,i}$ \label{line:lanczos_j=i-1_start}
\State $V_{*,i} \gets V_{*,i} - H_{i-1,i-1} V_{*,i-1}$ \label{line:lanczos_j=i-1_end}
\State $H_{i,i-1} \gets ||V_{*,i}||$
\State $V_{*,i} \gets \frac{V_{*,i}}{H_{i,i-1}}$
\EndFor
\State $H \gets H_{1:k, *}$ \Comment{discard extra row of $H$}
\State $V \gets V_{*, 1:k}$ \Comment{discard extra column of $V$}
  \end{algorithmic}
\end{algorithm}

We use Lemma~\ref{lm:error-bound} by incrementally increasing the number of Arnoldi iterations performed, $k$, until the
approximation error is smaller than a target accuracy.
The implemented algorithm starts with a small $k = 4$, and computes the corresponding error bound defined in Lemma~\ref{lm:error-bound}.
If the error bound satisfies the required accuracy, we use that value of $k$.
If not, we continue iterating, increasing $k$ by a factor of $1.1$ before the error is checked again.
The pseudocode is shown in Algorithm~\ref{alg:arnoldi_error}, with red lines indicating changes from the original Arnoldi process.
The \texttt{compute-error} function implements the error computation from Lemma~\ref{lm:error-bound}, which uses the $H$ matrix.

In our implementation, we target \old{an} \new{a} simulation accuracy of $10^{-6}$, the tolerance used in our LP solver.
In practice, the observed accuracy of the counter-examples produced upon reaching an unsafe state has been
significantly higher, as we will show in the evaluation.

\subsection{Krylov Simulations of Symmetric Matrices}
\label{ssec:lanczos}

A further improvement is possible when the system matrix is both sparse and \emph{symmetric}.
This may be the case when the dynamics matrix comes from a physical system due to the symmetry of many physical laws.
In this case, the Arnoldi iteration can be replaced by the more efficient Lanczos iteration~\cite{lanczos1950iteration,trefethen1997numerical}.
The difference between the two is that $H$ matrix in the symmetric case is both \emph{symmetric} and \emph{tridiagonal}.
This means that \old{that} step (2) in the algorithm, projecting out the previous orthonormal directions from the current vector,
only needs to be done for the previous two directions, and requires only a single dot product.
The Lanczos iteration is shown in Algorithm~\ref{alg:lanczos}.
Notice that the loop which projected out all the previous directions on lines~\ref{line:arnoldi_proj_start}-\ref{line:arnoldi_proj_end} of the original Arnoldi iteration in Algorithm~\ref{alg:arnoldi} is replaced by the $j=i-2$ case on lines~\ref{line:lanczos_j=i-2_start}-\ref{line:lanczos_j=i-2_end} and the $j=i-1$ case on lines~\ref{line:lanczos_j=i-1_start}-\ref{line:lanczos_j=i-1_end}.
This change reduces the computation time from $\mathcal{O}(k^2)$ to $\mathcal{O}(k)$.

Although the computation time is reduced with the Lanczos iteration, since the outputs are matrices of the same size as with Arnoldi, the memory required is basically the same as what was given in Equation~\ref{eq:memory_arnoldi}.
Some savings is possible if $H$ is stored as a sparse matrix, since $H$ is now tridiagonal rather than dense, but since typically $k \ll n$, this is savings is small.

\begin{algorithm}[t]
  \caption{Lanczos algorithm with projection and error control}\label{alg:lanczos_modified}
  \begin{algorithmic}[1]
\Require{normalized init $n \times 1$ vector $v$, $n \times n$ matrix $A$, \red{$o \times n$ projection matrix $C$, error target $\epsilon$}}
\Ensure{\red{$o \times k$ projected output matrix $\mathcal{P}=CV$}, $k \times k$ matrix $H$}
\State $V_{*,1} \gets v$
\State \red{$\mathcal{P}_{*,1} = C V_{*,1}$}
\State \red{$k \gets 4$}
\For{$i$ \textbf{from} $2$ \textbf{to} \red{$\infty$}}
\State $V_{*,i} \gets A ~ V_{*,i-1}$

\If{i > 2}
\State $H_{i-2,i-1} \gets H_{i-1,i-2}$
\State $V_{*,i} \gets V_{*,i} - H_{i-2,i-1} V_{*,i-2}$
\State \red{\texttt{free-memory}$(V_{*,i-2})$}
\EndIf
\State $H_{i-1,i-1} \gets (V_{*,i-1})^T V_{*,i}$
\State $V_{*,i} \gets V_{*,i} - H_{i-1,i-1} V_{*,i-1}$
\State $H_{i,i-1} \gets ||V_{*,i}||$
\State $V_{*,i} \gets \frac{V_{*,i}}{H_{i,i-1}}$
\State \red{$\mathcal{P}_{*,i} = C V_{*,i}$}
\RedIf{\red{$i = k$}} \red{\Comment{check error upon reaching $k$ iterations}}
\RedIf{\red{$\texttt{compute-error}(A, H_{1:i-1,*}) < \epsilon$}}
    \State \red{\textbf{break}}
    \RedElse
    \State \red{$k \gets \texttt{ceil}(1.1 * k)$}
\RedEndIf
\RedEndIf
\EndFor
\State $H \gets H_{1:k, *}$ \Comment{discard extra row of $H$}
\State \red{$\mathcal{P} \gets \mathcal{P}_{*, 1:k}$} \red{\Comment{discard extra column of $\mathcal{P}$}}
  \end{algorithmic}
\end{algorithm}

We propose a new modification to the Lanczos iteration that can save significant memory, when \old{its being} \new{it is} used for the verification problem.
Since we eventually project the result of $V_{k} e^{H_{k}t} \mathbf{e}_1$ onto the output space matrix $C$ (or the
transpose of the initial space matrix $E^T$), we propose to embed this projection inside the loop in the Lanczos algorithm, at each iteration.
The output of the iteration is then the much smaller $o \times k$ matrix $C V_k$ (or the $i \times k$ matrix $E^T V_k$).
This eliminates the need to store $V_{k}$, a potentially large $n \times k$ matrix, reducing the memory required for the algorithm.
The modified Lanczos iteration, which includes both this projection and the a posteriori error bound for selecting $k$, is shown in Algorithm~\ref{alg:lanczos_modified}.
As before, changes compared with the original Lanczos algorithm are in red.
With this improvement, the memory needed to perform the computation is reduced to:
\begin{equation}
\label{eq:memory_lanczos}
(3k + (n \times \min(i, o)) + 3n) \times \texttt{sizeof(double)}
\end{equation}
Importantly, compared with Equation~\ref{eq:memory_arnoldi}, there is no term where $k$ is multiplied by $n$.
This makes it possible to analyze high-dimensional systems with a much larger number of Krylov iterations, which may be needed for accuracy.
This will be needed in our evaluation when we evaluate a billion-dimensional system with $k=5932$ iterations.

%% file: scalability.tex
\subsection{Memory Scalability Limits}
\label{sec:scalability}

Several variables have been defined that impact the scalability of the proposed approach: 
$n$, the number of dimensions in the system dynamics, 
$i$, the initial space dimensions, 
$o$, the output space dimensions, 
$s$, the number of discrete time steps, and 
$k$, the dimension of the Krylov subspace used in the simulations, 
which is equal to the number of Arnoldi or Lanczos iterations needed.
Except for $k$, these are static variables, 
known before any computation is performed.
Using these variables, we can define the memory needed for the computation.

The amount of memory needed to store the basis matrix for all the steps is:
\begin{equation}
\label{eq:memory_basis}
o \times i \times s \times \texttt{sizeof(double)}
\end{equation}
%
%
Importantly, this limit is independent of the system dimensions $n$, which is why analysis with the proposed
approach can scale to extremely large systems.
In this case, even if one of $o = n$ or $i = n$, analysis may still be possible, as long as the product of 
$i$ and $o$ is manageable.

Next, if performing the Arnoldi iteration, we must also store $H_k$, a $k \times k$ matrix, and $V_k$, 
an $n \times k$ matrix.
The memory used by the Arnoldi algorithm is:
\begin{equation}
\label{eq:memory_arnoldi}
k \times (n + k) \times \texttt{sizeof(double)}
\end{equation}
The $k \times n$ factor in this equation is often the bottleneck, meaning that successful high-dimensional system
analysis cannot require a too high-dimensional Krylov subspace.

In the modified Lanczos iteration, $H_k$ is tridiagonal, and instead only the projection of $V_k$ is stored.
During the iteration, the current and previous two vectors of $V$ must be stored in order to be projected out, so a further 
factor of $3n$ is also needed.
The memory required for the Lanczos iteration is:
\begin{equation}
\label{eq:memory_lanczos}
(3k + (n \times o) + 3n) \times \texttt{sizeof(double)}
\end{equation}
If the transpose system simulation is used, $o$ is replaced by $i$.

Finally, the system matrix $A$ and initial space and output matrices $E$ and $C$ also need to fit into memory.
Even with a sparse representation, this can take non-negligible memory whens systems are very large.

%% file: evaluation.tex
\section{Evaluation}
\label{sec:evaluation}

We evaluate the proposed approach on several high dimensional
benchmarks\footnote{The source code and scripts to reproduce our evaluation results is available online: \url{http://stanleybak.com/papers/bak2019hscc_repeatability.zip}.}.
Measurements were performed using Amazon Web Services Elastic Computing Cloud (EC2), on a powerful
\texttt{m4.10xlarge} instance with 40 cores and a large 160 GB of memory that can be rented by the hour.
%
  %
  \new{Note that we perform analysis in discrete time, whereas \texttt{SpaceEx}~\cite{spaceex} and the decomposition method~\cite{bogomolov2018hscc}
    we compare against do dense time analysis.
  Although this requires a few more operations, specifically a bloating at the initial step using an approximation model, we do not expect
  qualitative runtime differences.} 

  \begin{figure*}[t]
    \centering
    \includegraphics[trim={0 0 0 0em}, clip, width=0.85\textwidth]{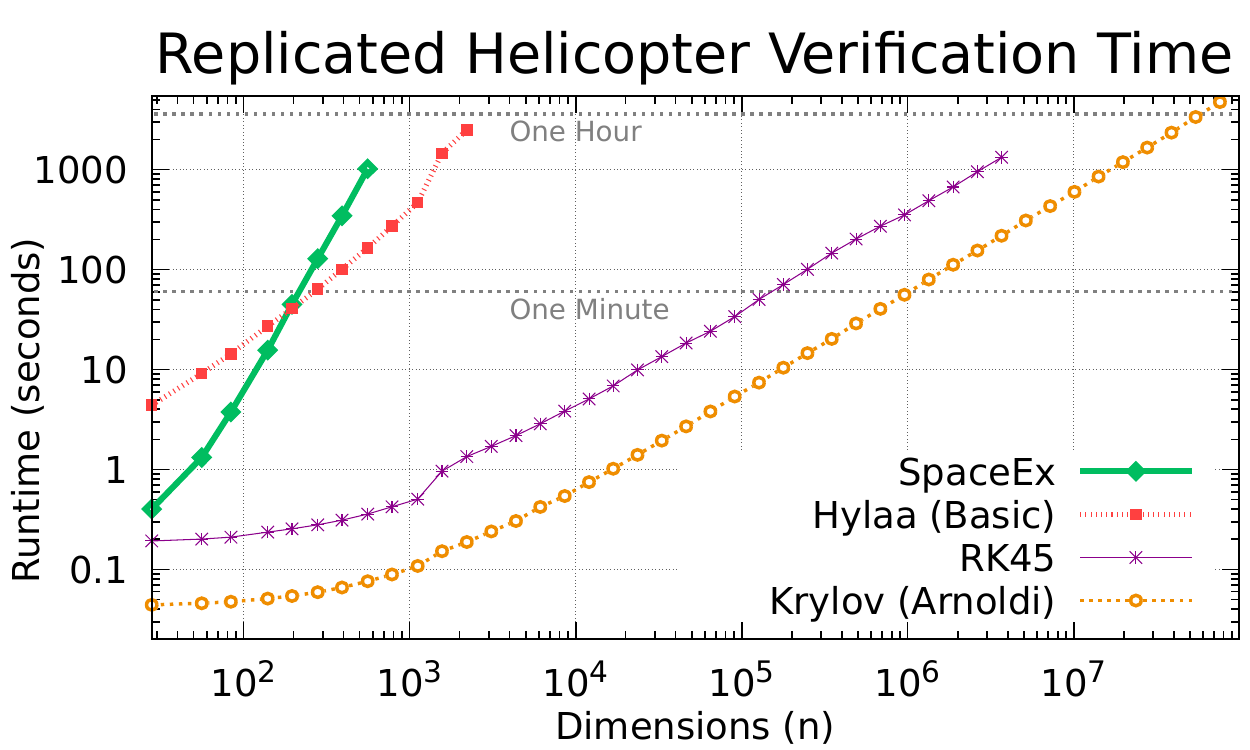}%
    \caption{The \texttt{Krylov} method scales better than the \texttt{RK45} numerical simulation method with input/output spaces,
    the full space method using numerical simulation used by \texttt{Hylaa}, and the space-time clustering
    scenario of the \texttt{SpaceEx} tool.}%
    \label{fig:heli_orig}%
\end{figure*}

\subsection{Modified Nodal Analysis (MNA5)}

We first verify a benchmark model based on a system from the field of electrical circuit analysis,
where the state variables relate to the node voltage and currents inside a
circuit~\cite{chahlaoui2002collection,tran2016large}.
Originally a DAE system, the dynamics matrix has been adapted to create a benchmark for verification using
ODE reachability methods.
As far as we are aware, this benchmark is the largest linear system ever verified~\cite{bak2017arch},
where full analysis of the safe version previously took a little over $24$ hours.
This model has also been investigated using a decomposition approach that uses a series of
two-dimensional projections to enable much faster analysis with modest overapproximation error~\cite{bogomolov2018hscc}.
Here, we apply the proposed algorithm which does not have overapproximation error and can provide
counter-examples when property violations are detected.

In this system, \new{the number of dimensions} $n = 10923$, \new{the number of output space dimensions} $o = 2$,
\new{the number of initial space dimensions} $i = 10$, and the number of steps is $20000$.
Our implementation selected \new{a Krylov subspace dimension} $k=63$ using the a posteriori error bound approach,
and verified the safe version of this system in $3.8$ \emph{seconds}.
The unsafe version of the benchmark was checked in $1.1$ seconds, with a counter-example at the same time in
the analysis as the earlier approach, at exactly step $1919$.

The initial state from the counter-example was then used to compute an external, high accuracy simulation of the system.
By comparing the final value in the external simulation versus the output variables assigned by the LP solver when the counter-example was found,
we can gauge our method's accuracy.
In this case, the relative error between the two was
$6.17 \times 10^{-9}$, demonstrating the accuracy of the proposed approach.

\subsection{Replicated Helicopter}
A tunable benchmark is created based on a 28-dimensional helicopter model and
controller originally released as an example system with the
SpaceEx tool\footnote{\url{http://spaceex.imag.fr/news/helicopter-example-posted-39}}.
The helicopter is copied multiple times within the same model, in order to create a verification problem that can
scale to an arbitrary number of dimensions.

In the replicated helicopter benchmark, the 28-dimensional helicopter model is copied $h$ times,
so that the number of dimensions $n = 28h$.
We take initial conditions from the \texttt{x8\_over\_time\_large} configuration,
where eight of the variables for each helicopter are initially intervals, making the dimension of
initial space $i = 8h$.
The error condition checks if the average of the $x_8$ variables is greater than $0.45$, $o = 1$.
Finally, the problem calls to verify up to time 30 with a step of 0.1, so that the number of steps is $300$.

Figure~\ref{fig:heli_orig} shows the scalability of the new approaches compared with
the \texttt{SpaceEx} tool and the basic approach implemented in the \texttt{Hylaa} tool.
We also tried to compare against the linear dynamics method in the Flow* tool~\cite{chen_rtss2014_flowstar},
but could not analyze the $h=1$ case due to the large uncertainty in the initial set.
Using the \texttt{stc} scenario~\cite{frehse2013flowpipe} of \texttt{SpaceEx}~\cite{spaceex}, the largest system
successfully analyzed had $h=20$ (560 dimensions) and took 17 minutes (larger systems crashed).
The basic approach implemented in the \texttt{Hylaa} tool, which computes the full $n \times n$ basis matrix using numerical
simulations, verified a system $h=79$ (2212 dimensions) in 42 minutes (larger systems
had a memory error).
Using the proposed input / output spaces with Runge-Kutta numerical simulations, the \texttt{RK45} method scaled
up to $h=131389$ (3.6 million dimensions) in about 22 minutes (larger systems had a memory error).
Finally, combining both input / output spaces and Krylov subspace simulations, the \texttt{Krylov} approach analyzed
the system with $h=2714654$ (76 million dimensions) in 79 minutes, without memory errors.
Since the initial space dimension $i$ grows as the number of helicopters $h$ increases, the LP solving step takes
increasing amounts of time for this benchmark.
This explains why the slopes for \texttt{RK45} and \texttt{Krylov} are similar: the LP solving step
has become the bottleneck.

%
%
%
%

\subsection{Symmetric 3D Heat Diffusion}
\label{ssec:heat3d}
The third benchmark considered is a 3D Heat Diffusion system taken from the field of
partial differential equations (PDEs).
This benchmark is based on a 2D version that has previously been analyzed up to a
$50 \times 50$ mesh ($2500$ dimensions)~\cite{han2005formal,han2006reachability}.
The problem is to examine the temperature at the center point of a $1.0 \times 1.0 \times 1.0$ block, where
one edge of the block is initially heated.
As before, all of the sides of the block are insulated except the $x=1.0$ edge, which allows
for heat exchange with the ambient environment with a heat exchange constant of $0.5$.
%
A heated initial region is present in the region where $x \in [0.0, 0.4]$, $y \in [0.0, 0.2]$, and $z \in [0.0, 0.1]$.
The heated region temperature is between $0.9$ and $1.1$, with the rest of material
initially at temperature $0.0$.
The system dynamics is given by the heat equation PDE $u_t = \alpha^2(u_{xx} + u_{yy} + u_{zz})$,
where $\alpha = 0.01$ is the diffusivity of the material, as in the previous work.

\begin{figure}[tb]%
  \centering
\includegraphics[width=0.65\textwidth]{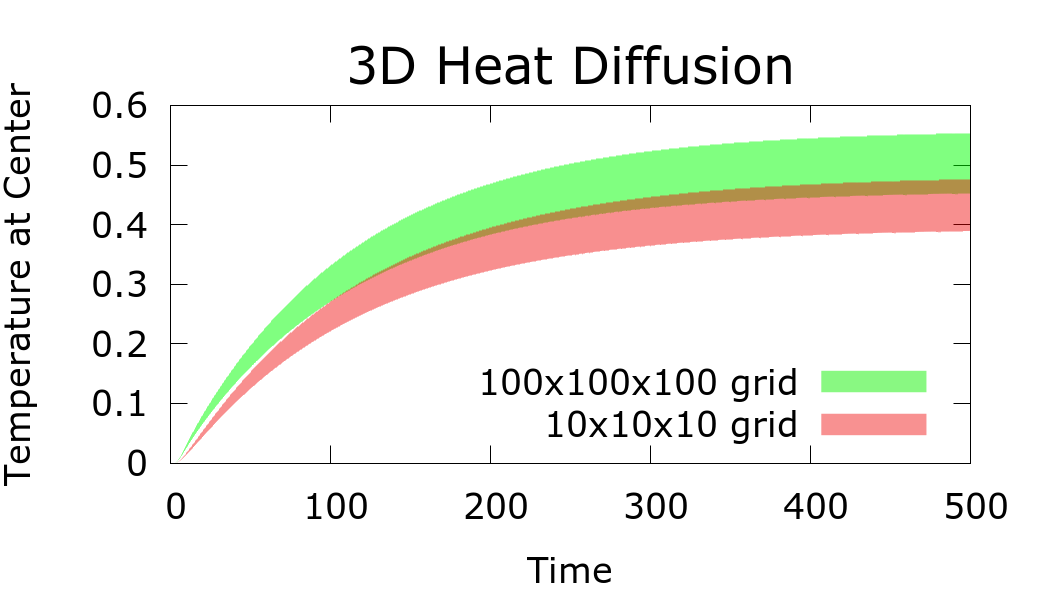}%
\caption{The maximum temperature at the center point occurs around time 15.}%
\label{fig:heat_sym_reach}%
\end{figure}

A linear state space model of the system is obtained using the semi-finite difference method~\cite{farlow1993partial},
discretizing the block with an $m \times m \times m$ grid.
This results in an $m^3$-dimensional linear system describing the evolution of the temperature at each mesh point.

Due to the initially heated region, we expect the temperature at the center of the block to first increase, and then
decrease due to the heat loss along the $x=1$ edge.
Further, there may be error due to the space discretization step, so if $m$ is too small,
the model does not accurately predict the behavior of the PDE.
We can see both of these effects by computing and plotting the reachable states,
as shown in Figure~\ref{fig:heat_sym_reach}.
Since the peak temperature happens at around time $t=15$, we perform further analysis by running the system
with max time $T=20.0$ and step size $\delta=0.02$, making the number of steps $1000$.
This system presents a particularly good case for our analysis method, since $i=1$, $o=1$, and the dynamics
matrix is symmetric which allows us to use the Lanczos iteration.

The runtimes and temperatures reachable for various values of $m$ are given
in Table~\ref{tab:heat_sym}.
Accurate analysis requires high dimensions, motivating the need for
the types of analysis methods developed in this paper.
The $1000 \times 1000 \times 1000$ version can be analyzed using our approach in about 30 hours of
computation time.
Over 95\% of the runtime was spent in the Lanczos iteration,
indicating that we optimized the correct operation.
In this case, each of the billion rows of the $A$ matrix generally has $7$ entries, so that simply storing the
elements of the matrix ($8$ bytes per double-precision number) consumes $56$ GB of RAM.
%
Further, since a 5932-dimensional Krylov subspace is needed for sufficient numerical accuracy,
the unmodified \old{Krylov}\new{Lanczos} iteration would be infeasible for this system,
as it would require storing 5932 vectors for the $V$ matrix, each of which contains a billion numbers ($8$ GB each),
for a total memory requirement of 46 TB (recall Equation~\ref{eq:memory_arnoldi}).

\begin{table}[t]%
  \centering
  \caption{3D Heat Diffusion with $n=m^3$ Dimensions}%
  \label{tab:heat_sym}%
\begin{tabular}{@{}lllll@{}}%
\toprule%
\bf{$m$} ~~~ & $T_{\textnormal{max}}$ ~~~ & $k$ ~~~~ & \bf{Lanczos} & \bf{Arnoldi} \\%
\midrule%
\input{figures/heat_sym_summary.dat}\bottomrule%
\end{tabular}%
\end{table}

\begin{figure}[tb]
  \centering
  \includegraphics[width=0.75\textwidth]{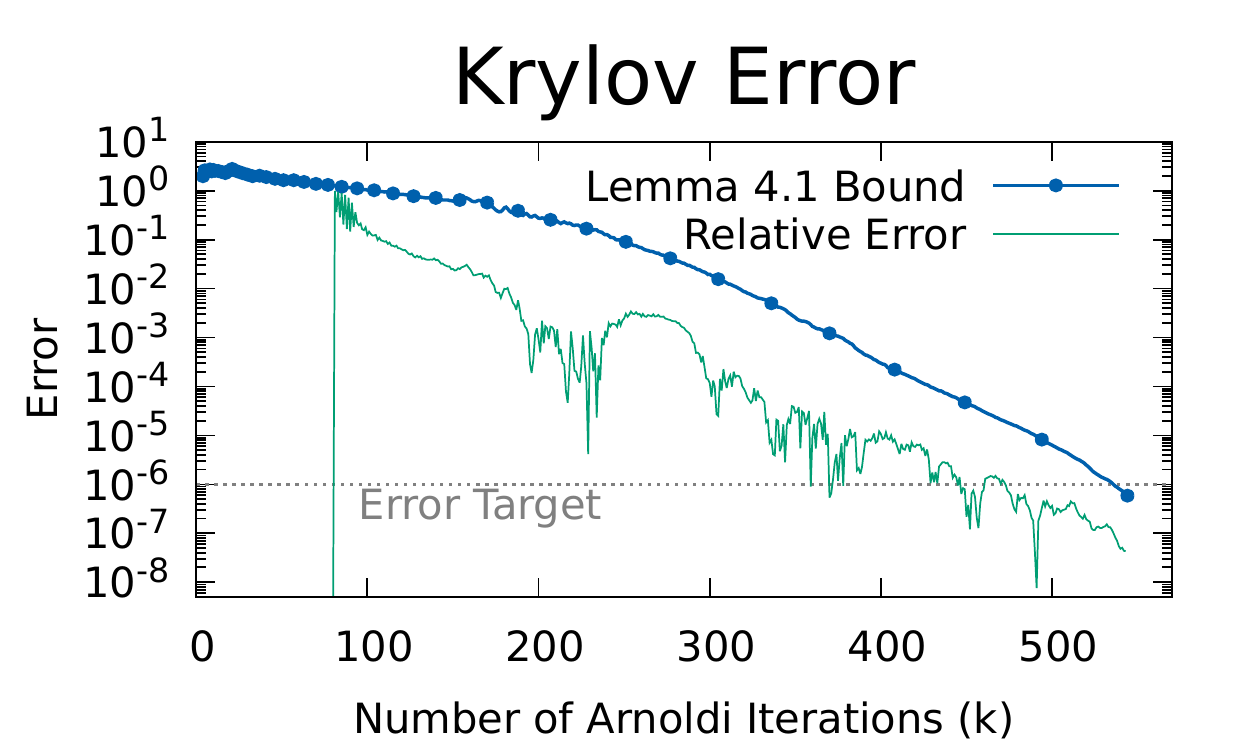}
\caption{A 544-dimensional Krylov subspace exceeds the desired error target of $10^{-6}$ for the 100x100x100 \new{(one million dimensional)} 3D Heat diffusion system.}%
\label{fig:krylov_error}%
\end{figure}

Lastly, we examine the error bound from Lemma~\ref{lm:error-bound} for the 100x100x100 version of this benchmark,
as the dimension of the Krylov subspace $k$ is increased.
While performing the Arnoldi or Lanczos algorithm, our implementation periodically checks the current error.
Once $k=544$ iterations have been performed, the computed error bound is $5.8 * 10^{-7}$, which is below the desired error
threshold of $10^{-6}$.
The plot is shown in Figure~\ref{fig:krylov_error}.
\new{The blue line is the error bound computed using Lemma~\ref{lm:error-bound} at each iteration, and the points indicate where the bound gets sampled.
  The thinner green line is the relative error of the projected simulation between iterations $k$ and $k+1$, which
  provides an error \emph{estimate} that was used as a stopping criteria in previous work~\cite{bak18arch}.
  Notice that the old error bound can not be used when the number of iterations is low, as the projected simulation onto the output variables is zero
  when $k$ is small.
  For this system, using the old bound might reach the $10^{-6}$ threshold earlier and terminate prematurely,
  although the number of iterations for both bounds is similar, within a factor of two.
}
Recall from the discussion after Equation~\ref{eq:error-bound} that the
a priori bound was unusable for this system, even with $k=10^6$.

%
%

%% file: figures/heat_sym_summary.dat
10 & 0.02966 & 63 & 0.5s & 0.3s \\
20 & 0.01716 & 94 & 0.5s & 0.4s \\
50 & 0.01161 & 211 & 1.3s & 6.7s \\
100 & 0.01005 & 474 & 7.4s & 6m28s \\
200 & 0.00933 & 711 & 1m40s & 2h3m \\
500 & 0.00891 & 1599 & 45m55s & - \\
1000 & 0.00877 & 3597 & 13h26m & - \\

%% file: related_work.tex
\section{Related Work}
\label{sec:related_work}
The proposed method uses convergent numerical schemes to compute simulations as part of a
verification procedure.
Convergent numerical schemes have been used before to approximate reachable sets of nonlinear hybrid systems,
in particular, level-set methods that approximate solutions to Hamilton-Jacobi PDEs~\cite{mitchell2000level,tomlin2003,bansal2017hamilton}.
These methods compute reachable states with a grid over the state space, and
in the limit at the number of grid points increases, the computed result approaches the true solution.

Other methods for this class of systems have used simulations for formal analysis,
where individual executions are bloated according to model-specific discrepancy
functions~\cite{fan2015bounded}, as implemented in tools such as C2E2~\cite{duggirala2015c2e2,fan2016automatic}.
Another analysis approach for nonlinear systems uses Taylor models,
such as those in Flow*~\cite{chen_rtss2014_flowstar}, which
can scale to around ten real variables~\cite{chen2015benchmark}.
For affine systems, as recently as 2011 the state-of-the-art for
reachability computation was on the order of a hundred real variables~\cite{spaceex}.

Our work uses the Krylov subspace to simulate high-dimensional systems, which is often also used in
model order reduction methods~\cite{Antoulas01asurvey}.
Notice that in our case, since each simulation has a different Krylov subspace,
\emph{there is no single reduced order model that can be constructed and analyzed}
(we are not creating a low-dimensional abstraction of the system).
Model-order reduction approaches verify a smaller dimensional model~\cite{chou2017study},
and can sometimes use an error bound to compute a guaranteed overapproximation of
the original full-order system~\cite{han2004reachability,han2005formal,tran2017order}.
Such approximation methods may be formalized as sound abstractions or developed in the context of
approximate simulation and bisimulation relations~\cite{girard2007approximation,girard2011approximate}.
Model order reduction methods have verified linear systems with on the order of
a thousand real variables.

Our approach builds on the basic verification approach used in the Hylaa
tool~\cite{bak2017hscc}, which has verified systems with up to ten thousand
dimensions~\cite{bak2017arch,bak2017cav,bak18arch}.
We scale to larger systems here by leveraging initial and output spaces and using Krylov subspace methods for numerical simulation.
The basic approach here is also related to the symbolic orthogonal projection
method~\cite{hagemann2014reachability}, where the
current-time variables in our approach could be considered the variables onto which we are computing the projection.

Recent work on reachablity with Krylov methods~\cite{althoff2017reachability} has used better a priori bounds to
reduce the number of required Arnoldi iterations, compared with earlier work~\cite{han2006reachability}.
Further, more efficient methods exist which perform Krylov simulations in multiple stages~\cite{krylov_algorithm_variable_step},
rather than only from the initial time as in this work.
Integrating these into our approach could further improve our result by reducing the required dimension of the Krylov subspace.

%
%

%
%
%
%
For large systems, the only analysis option we are aware of beyond simple simulation would be falsification
methods~\cite{annpureddy2011s,donze2010breach,rizk2008continuous}, which run individual simulations trying to
optimize towards an unsafe region.
Unlike our approach, these methods do not exhaustively explore the state space.


%% file: conclusion.tex
\section{Conclusion}
\label{sec:conclusion}

The state-space explosion problem usually prevents analysis of high-dimensional affine systems.
In order to achieve scalability, we exploited up to four types of problem structure:
(i) a small dimension of initial states, (ii) a small dimension of the output space,
(iii) the sparsity of the system $A$ matrix, and (iv), optionally, the symmetry of the $A$ matrix.
When problems have this structure, we have shown it is possible to efficiently perform verification or 
plot projections of the reachable states despite a large number of state variables.
As the structure assumptions are violated, the proposed approach degrades gracefully, requiring more 
computation time and memory depending on the degree of the violation.
We have evaluated our approach on several large benchmarks,
including a 3D Heat Diffusion system with one billion continuous state variables.
Prior to this work, no existing method for affine systems has demonstrated scalability beyond 
a few thousand variables.